\documentstyle{amsppt}
\input epsf

\magnification=\magstep1
\vsize=9 true in
\hsize=6.5 true in
\topmatter
\title
Transition Mean Values of Real Characters
\endtitle
\author
J.B.~Conrey \\
D.W.~Farmer \\
K.~Soundararajan
\endauthor
\thanks Research of the first author supported in
part by a grant from the NSF.  
Research of all three authors supported in 
part by the American Institute of Mathematics.
\endthanks
\address
 American Institute of Mathematics, 
Palo Alto, CA 94306
\endaddress
\address
School of Mathematics,  Institute for Advanced Study,
Princeton, NJ 08540
\endaddress           
\abstract
We evaluate the real character sum $\sum_m\sum_n \fracwithdelims() mn$ 
where the two sums are of approximately the same length.   
The answer is surprising.
\endabstract
\endtopmatter
\document

\NoBlackBoxes

\def\hb{\relax}
\def\mod{\bmod}

\def\Q {{\Bbb Q}}

\def\intl{\int\limits}

\def\vp{\varphi}
\def\al{\alpha}

\def\({\left(}
\def\){\right)} 
\def\js{\fracwithdelims()}

\head
1.  Introduction and statement of results
\endhead

\noindent Let
$$
S(X,Y):=\sum_{{ m\le X}\atop m\text{  odd} }
\sum_{{n\le Y}\atop n\text{  odd}}
\js mn ,
$$
where $\fracwithdelims() mn$ is the Jacobi symbol.  
Our goal is to obtain an asymptotic formula
for $S(X,Y)$. 
We will see that this is straightforward except
when $X$ and $Y$ are of comparable size.

First we give asymptotic formulas valid
for $Y=o(X/\log^{}X)$ or
$X=o( Y/\log^{}Y)$.
An easy application
of the P\'olya--Vinagradov inequality shows that
$$
\sum_{{m\le X}\atop m \hb\text{  odd}} \js mn
=
\cases
{\displaystyle{X\over 2}{\varphi(n)\over n}}+ O(X^\varepsilon) &
\hbox{if } n =\square \\
O(n^{1\over 2} \log n) & \hbox{if } n \not=\square \\
\endcases
$$
and
$$
\sum_{{n\le Y}\atop n \hb\text{ odd}} \js mn
=
\cases
{\displaystyle{Y\over 2}{\varphi(m)\over m}}+ O(Y^\varepsilon) 
&\hbox{if } m =\square \\
O(m^{1\over 2} \log m) & \hbox{if } m \not=\square ,\\
\endcases 
$$
where $\square$ represents the square of
a rational integer.

It follows that
$$
\eqalignno{
S(X,Y)&=
\sum_{{ n\le Y}\atop n = \text{  odd }\square}
\({X\over 2}{\varphi(n)\over n} +O(X^\varepsilon)\)
+O(Y^{3\over 2}\log Y)\cr
&={2\over \pi^2} XY^{1\over 2} + 
O(Y^{3\over 2} \log Y + Y^{1\over 2}X^\varepsilon+X\log Y), &(1.1)\cr
}
$$
and similarly,
$$
S(X,Y)={2\over \pi^2} X^{1\over 2} Y + 
O(X^{3\over 2} \log X + X^{1\over 2}Y^\varepsilon + Y\log X).
\eqno(1.2)
$$
Equation (1.1) provides an asymptotic formula for 
$S(X,Y)$ when $Y=o(X/\log X)$, and (1.2) when $X=o( Y/\log^{}Y)$.
The range when $X$ and $Y$ are of comparable size  marks a 
transition in the behavior of $S(X,Y)$, and our 
object here is to understand this
transitory phase.

\proclaim{Theorem 1}  Uniformly for all large $X$ and $Y$, we have
$$
S(X,Y) ={2\over \pi^2} C\(\frac{Y}{X}\) X^{3\over 2} + O\((XY^{7\over 16}+
YX^{7\over 16})\log XY\),
$$
where for $\alpha\ge 0$ we define 
$$
C(\alpha) = \sqrt{\alpha}+ \frac{1}{2\pi} \sum_{k=1}^{\infty} 
\frac{1}{k^2} \intl_0^{\alpha} \sqrt{y} 
\biggl(1-\cos \(\frac{2\pi k^2}{y}
\)+\sin\(\frac{2\pi k^2}{y}\)\biggr) dy.
$$
An alternate expression for $C(\alpha)$ is 
$$
C(\alpha)={\alpha}+ \alpha^{\frac 32}\frac{2}{\pi} \sum_{k=1}^{\infty} 
\frac{1}{k^2} \intl_0^{1/\alpha} 
\sqrt{y} \sin\(\frac{\pi k^2}{2y}\) dy.
$$
\endproclaim

To assist the reader in understanding the function $C(\alpha)$,
graphs of $C(\alpha)$ and $C'(\alpha)$ are presented in Section 6.  

The first expression for $C(\alpha)$ shows, upon integrating by parts, that 
$$
C(\alpha) =\sqrt{\alpha} + {\pi \over 18}\alpha^{3\over 2}
+O(\alpha^{5\over 2}) \qquad \text{as  } \al \to 0.
$$
Similarly, the second expression for $C(\al)$ gives the limiting behavior
$$
C(\al) = \al+ O(\al^{-1}) \qquad \text{as  } \al \to \infty.
$$
Note that in these limiting cases, the value of 
$C(\al)$ approaches that given by the $n=\square$ terms (as $\al\to 0$) 
and the $m=\square$ terms (as $\al\to\infty$).  From these limiting behaviors 
(or (1.1) and (1.2)) we see that $C(Y/X)X^{3\over 2}$ is of size 
$XY^{1\over 2} + YX^{1\over 2}$, so that the error term of 
Theorem~1 is always smaller than the main term.  We shall leave to the 
reader the problem of showing that our two expressions for $C(\al)$ agree:
this is an exercise in the Poisson summation formula.

If quadratic reciprocity said $\js mn = \js nm$ for all $m$ and~$n$,
then we would have a functional equation 
$C(\alpha)=\alpha^{3\over 2} C(1/\alpha)$.  
Plainly this does not hold,
although the above expressions do show a
relationship between $C(\alpha)$ and $\alpha^{3\over 2}C(1/\alpha)$.
If the sum defining $S(X,Y)$ had been restricted to
$m\equiv n\equiv 1\pmod 4$, then such a functional 
equation would hold.
 
Differentiating our second expression for $C(\al)$ term by term we 
obtain
$$
\(\alpha^{-\frac 32}C(\al)\)' =    
-\frac{1}{2} \al^{-\frac 32} -
\al^{-\frac{5}{2}}
f\(\frac{\mathstrut \alpha}{2} \),
$$
where
$$      
f(x)=\frac 1{\pi}
\sum_{k\not=0}
\frac{1}{k^2}                                        
\sin\( \pi k^2x \)
.          
$$
The function $f(x)$ is commonly called 
``Riemann's nondifferentiable function,''
so named because Weierstrass reported that Riemann suggested
it as an example of a continuous function which is not differentiable.
A considerable amount of work has been done investigating the 
differentiability properties of $f(x)$.
Hardy~\cite{H} showed that it is not differentiable at $x=s$
if $s$ is irrational or if $s=\frac pq$ with $p$ or $q$ even.
Gerver~\cite{G} gave a long elementary proof that 
$f'(p/q)=-1$ if $p$ and $q$ are odd, so Riemann's
assertion is not quite correct.
For an interesting survey on Riemann's function,
see Duistermaat~\cite{Du}.  In Section~6 we show that
determining  the differentiability of $f(x)$ at a
rational point is a straightforward exercise in the 
Poisson summation formula; our approach appears to be similar
to that of Smith~\cite{Sm}.
It seem surprising that the asymptotics of a natural object
like $S(X,Y)$ should involve non--smooth functions!

In the following section we explain our motivation for studying
the sum $S(X,Y)$, and we give a generalization of Theorem~1.  
In Section~3 we do some preliminary reductions
and identify the main terms and error terms in the sum.  These are
evaluated and estimated in Sections~4 and~5.
In Section~6 we present graphs of~$C(\al)$ and~$C'(\alpha)$,
and we determine at which rationals $C'(\alpha)$ is differentiable.

\head
2.  Motivation: mollifying $L({1\over 2},\chi_{d})$
\endhead

\noindent The motivation for studying $S(X,Y)$ came 
from the third author's proof \cite{S}
that $L({1\over 2},\chi_{d})\not=0$ 
for a positive proportion of fundamental discriminants $d$,
where $\chi_{d}$ is the real primitive character to the modulus $d$.
Jutila~\cite{J} showed that there exist positive constants $c_1$ and $c_2$ such that 
$$
\leqalignno{
\sum_{|d|<X} L(\textstyle{1\over 2},\chi_{d}) &\sim c_1 X\log X\cr
&&\hbox{and}\cr
\sum_{|d|<X} |L(\textstyle{1\over 2},\chi_{d})|^2 &\sim c_2 X\log^3 X ,
}
$$
where both sums range over fundamental discriminants.  It follows
from the above formulas and Cauchy's inequality that the
number of $|d|<X$ with $L({1\over 2},\chi_{d})\not=0$ exceeds $X/\log X$.

The approach used to obtain the nonvanishing of
$L(\textstyle{1\over 2},\chi_{d})$ for a positive proportion of $d$
was to consider a ``mollified'' sum of $L(\textstyle{1\over 2},\chi_{8d})$.
Let
$$
M(d)=\sum_{\ell\le M}\lambda(\ell)\sqrt{\ell}  \js{8d}\ell ,
$$
where $\lambda(\ell)$ is chosen so that
$$
S_1:=\sum_{|8d|<X} L({\textstyle{1\over 2}},\chi_{8d}) M(d)
\ \ \ \ \ \ \ \ \ \ 
\hbox{and}
\ \ \ \ \ \ \ \ \ \
S_2:=\sum_{|8d|<X} |L({\textstyle{1\over 2}},\chi_{8d})M(d)|^2 
$$
are both of size $X$.  By Cauchy's inequality this implies that
$L(\textstyle{1\over 2},\chi_{8d})\not=0$ for a positive proportion of~$d$.
The optimal choice for $\lambda(\ell)$ is determined in~\cite{S}.
The answer is complicated, so suffice it to say that
$\lambda(\ell)$ is supported on the odd integers, where
$$
\lambda(\ell)
\ \ \ \ 
\hbox{ is roughly proportional to }
\ \ \ \ 
\frac{\mu(\ell)}\ell 
\frac{\log^2(M/\ell)}{\log^2 M}
\frac{\log(X^{3\over 2}M^2\ell)}{\log M}.
$$ 
This leads to the result that 
$L(\textstyle{1\over 2},\chi_{8d})\not=0$
for at least $\frac 78 $ of all squarefree integers~$d$.

The most difficult part of the above argument is the evaluation
of a certain ``off--diagonal'' contribution to the main term.
This involves finding an asymptotic formula for an expression
of the form
$$
\Sigma_\ell(X):=\sum_{m,n}
F(m,n,\ell,X)
\js {m}{\ell n},
$$
for some explicit function~$F$.  In this paper we supress the
function~$F$, and we find that the resulting sum retains
the interesting features of the corresponding sum considered 
in~\cite{S}.

Motivated by the sum $\Sigma_\ell(X)$ we also
consider the
slightly
more general sum
$$
S_\ell(X,Y):=
\sum_{{ m\le  X}\atop m\hb{\ odd}}
\sum_{{ n\le Y}\atop n\hb{\ odd}}
\js m{\ell n} .
$$

\proclaim{Theorem 2}  Let $\ell$ be an odd squarefree integer.
There exists a constant $C_\ell(\alpha)$ such that, uniformly for all large $X$ and $Y$,
$$
S_\ell(X,Y)\sim
{1\over \sigma(\ell)}
{2\over\pi^2}
C_\ell\(\frac YX\) 
X^{3\over 2},
$$
where $\sigma(\ell)$ is the divisor sum function.
We have 
$
C_\ell(\alpha) =
C(\alpha\ell) 
$,
where $C(\alpha)$ is given in Theorem~1.
\endproclaim

The proof of Theorem 2 will be omitted because it closely follows
the proof of Theorem~1.

\head
3. Initial reductions
\endhead

\noindent When $Y\le X^{16 \over 17}$ or $X\le Y^{16\over 17}$, 
Theorem 1 follows immediately from (1.1) and (1.2) and the limiting 
behaviors of $C(\al)$.  We assume below that 
$Y^{16\over 17} 
\le X\le Y^{17 \over 16}$.  
In place of $S(X,Y)$ it is technically easier
to consider the smoothed sum
$$
\Cal S (X,Y):=\sum_{m \text{  odd}} \ \ \sum_{n\text{  odd}} 
    \js mn H\biggl(\frac{m}{X}\biggr) \Phi\biggl(\frac{n}{Y} \biggr).
$$
Here $H$ and $\Phi$ are smooth functions supported
in $(0,1)$, 
satisfying $H(t)=\Phi(t)=1$ for $t\in (1/U,1-1/U)$, and such that 
$H^{(j)}(t),\ \Phi^{(j)}(t) \ll_j U^j$ for all integers $j\ge 0$.  
The parameter $U$ will later be chosen to equal 
$(XY)^{\frac 25}(X+Y)^{-\frac 35}$.  

Using the P{\'o}lya--Vinogradov inequality in a way similar to the
argument described in the Introduction, it is easy to see that 
$$
|S(X,Y) -{\Cal S}(X,Y)| 
\ll \frac{X^{3\over 2} +Y^{3\over 2}}{U}\log XY.
 \eqno(3.1)
$$
With our choice of $U$ this is seen to be smaller than the error term.

We evaluate ${\Cal S}(X,Y)$ by applying the Poisson summation formula to 
the sum over $m$ in ${\Cal S}(X,Y)$.  For all odd integers $n$ and 
all integers $k$, we  introduce the Gauss--type sums 
$$
\tau_k(n):=\sum_{a\pmod n} \js an e\biggl(\frac{ak}{n}\biggr)
=:\({1+i\over 2} + \js {-1}n {{1-i\over 2}}\) G_k(n),
$$
where $e(x):=e^{2\pi i x}$ as usual. 
We quote Lemma 2.3 of [S] which determines $G_k(n)$.  

\proclaim{Lemma 1} If $(m,n)=1$ then $G_k(mn)=G_k(m)G_k(n)$.  Suppose
that $p^a$ is the largest power of $p$ dividing $k$ 
(put $a=\infty$ if $k=0$).  Then for $b\ge 1$ we have
$$
G_k(p^b)=
\cases
0&\text{if $b\le a$ is odd}\cr
\vp(p^b) & \text{if $b\le a$ is even}\cr
-p^a&\text{if $b=a+1$ is even}\cr
{\js{k/p^a}{p}} p^a\sqrt{p} & \text{if $b=a+1$ is odd}
\cr
0&\text{if $b\ge a+2$}.\cr
\endcases
$$
\endproclaim

By Poisson summation we have (see section 2.4 of \cite{S} for
details):
$$
\sum_{m\text{   odd}} \ \  \js mn H\({m\over X}\)
={X\over 2 n}\js 2n \sum_{k=-\infty}^\infty 
(-1)^k \tau_k(n) \widehat H\({k X\over 2 n}\).
$$
Expressing $\tau_k$ in terms of $G_k$, using the relation 
$G_k(n)=\js {-1}n G_{-k}(n)$, 
and recombining the $k$ and $-k$ terms we may rewrite the above as
$$
{X\over 2 n}\js 2n \sum_{k=-\infty}^\infty (-1)^k G_k(n)
\,\widetilde{H}\({kX\over 2n}\),
$$
where
$$
\widetilde{H}(\xi):=
{1+i\over 2} \widehat H (\xi)+{1-i\over 2} \widehat H(-\xi).
$$  

These manipulations show that 
$$\eqalignno{
{\Cal S}(X,Y) &= \frac{X}{2} \sum_{k=-\infty}^{\infty} 
\sum_{n \text{  odd}} \ \ (-1)^k  \fracwithdelims(){2}{n}
\frac{G_k(n)}{n} \Phi\biggl(\frac{n}{Y}\biggr)
\widetilde{H} \biggl(\frac{kX}{2n}\biggr) \cr
&= M+R,&{(3.2)}
}$$
where $M$ isolates the terms when $k=2\square$, and $R$ includes the 
remaining terms.  That is,
$$\eqalign{
M:&= {X\over 2}\sum_{k=0}^\infty \,\sum_{n\text{  odd}}  \ \ 
{G_{2k^2}(n)\over n}\js 2n \Phi\({n \over Y}\) 
\,\widetilde{H}\({k^2X\over n}\)\cr
&={X\over 2}\sum_{k=0}^\infty M_k,
}$$
say, and
$$
R:={X\over 2}\sum_{k=-\infty \atop k\not=2\square}^\infty(-1)^k 
\sum_{n\text{ odd}} \ \   {G_k(n)\over n} \js 2n \Phi\({n\over Y}\)
\,\widetilde{H}\({kX\over 2n}\).
$$
We will see that $M$ gives the main term and $R$ is an error term.

\head
4. The remainder term $R$
\endhead

We will require some simple estimates on $\widehat H(t)$
and $\widetilde H(t)$.  These follow by integration by parts
and our assumptions on $H$ and $\Phi$.  We have
$$
|\widehat H(t)|,\ \ |\widetilde H(t)|,\ \ |(\widetilde H(t))^{\prime}|
\ll_j U^{j-1} |t|^{-j}
\eqno(4.1)
$$
for all integers $j\ge 1$, and all real $t$, and
$$
{\widetilde H}(\xi) =  \frac{1-\cos (2\pi \xi)+\sin(2\pi \xi)}{2\pi \xi}
+O\left( \frac{1}{U} \right). 
\eqno(4.2)
$$

We handle the remainder term $R$ using the following Lemma which 
exhibits cancellation in the sum $G_k(n) \js 2n$ when $2k\neq \square$. 

\proclaim{Lemma 2}  If $k\neq 2\square$ then 
$$
\sum_{n\le x \atop n\text{ odd}} \frac{G_k(n)}{\sqrt{n}} \js 2n
\ll |k|^{1\over 4} \log (2|k|) d(k^2) x^{1\over 2},
$$
where $d(k^2)$ is the number of divisors of $k^2$.
\endproclaim

Before proving the Lemma we note the bound it gives for $R$.  By partial 
summation and Lemma~2 we have
$$
\align
\sum_{n\text{ odd}} \frac{G_k(n)}{n} \js 2n \Phi\(\frac nY\)
\widetilde H\(\frac{kX}{2n}\) 
\ll |k|^{1 \over 4} & \log (2|k|) d(k^2) \\
& \times\intl_0^Y \sqrt{t} 
\biggl|\biggl(\frac{1}{\sqrt {t}} \Phi\(\frac tY\) 
\widetilde H\(\frac{kX}{2t}\)\biggr)^{\prime}\biggr|dt,\\
\endalign
$$
and using (4.1) with $j=3$ this is 
$$
\ll |k|^{1\over 4} \log (2|k|) d(k^2) \frac{U^2 Y^2}{k^2 X^2}.
$$
Summing over all $k\neq 2\square$ we obtain
$$
R\ll \frac{U^2 Y^2}{X}. 
\eqno(4.3)
$$

\demo{Proof of Lemma 2} We write $n=rs$ where $r$ and $s$ are odd with $s$ 
coprime to $k$ and $r$ divisible only by primes dividing $k$.  By Lemma
1, $G_k(n) =G_k(r)G_k(s)= G_k(r) \mu^2(s)\sqrt{s}\js ks$.  Further 
note that $G_k(r)=0$ unless $r|k^2$, and at any rate $|G_k(r)|\le r$.  Thus 
our desired sum is 
$$
\ll \sum\Sb r|k^2 \\ r\le x\endSb \sqrt{r} \biggl|\sum\Sb s\le x/r 
\endSb \mu^2(s)\js {2k}{s} \biggr|.
$$
Expressing $\mu^2(s)=\sum_{d^2|s}\mu(d)$, and using the P{\' o}lya-Vinogradov
inequality (since $\js {2k}{\cdot}$ is a non-principal character 
with conductor $\le 8|k|$) we obtain
$$
\align
\sum\Sb s\le x/r \endSb \mu^2(s) \js{2k}{s} &\ll 
\sum_{d\le \sqrt{x/r}} \biggl|\sum_{t\le x/rd^2} \js{2k}{t}\biggr| 
\ll \sum_{d\le \sqrt{x/r}} \min\biggl(\frac{x}{rd^2}, |k|^{1\over 2} 
\log (2|k|)\biggr) \\
&\ll |k|^{1\over 4} (\log |2k|) x^{1\over 2} r^{-{1\over 2}}.\\
\endalign
$$
Using this in our previous display we obtain Lemma~2.
\enddemo
 
\head
5. The main term $M$
\endhead

\noindent First consider the case $k=0$.  It follows straight
from the definition that $G_0(n)=\vp(n)$ if $n=\square$ and
$G_0(n)=0$ otherwise.  Thus 
$$
M_0 ={\widehat H}(0) \sum_{n=\text{ odd } \square}\ \ 
{\vp(n)\over n} \Phi\({n \over {Y}}\)
=   \widehat H(0) \frac{2\sqrt{Y}}{\pi^2}  
 \intl_0^\infty {\Phi(y)\over \sqrt{ y}}dy  
+ O(\log Y).
$$
The second step follows from 
$$
\sum\Sb n \le x\\ n= \text{ odd }\square\endSb \ \ \frac{\vp(n)}{n} 
= \frac{4}{\pi^2} \sqrt{x} + O(\log x),
$$ 
and then using partial summation.

Now suppose $k\ge 1$.  Since $n$ is odd, by changing variables
in the sum defining $G_k$ we have $G_{2k}(n)\js 2n = G_k(n)$.
This gives
$$
M_k=\sum_{n\text{ odd}} \ \  \frac{G_{k^2}(n)}{n}
\Phi\({n \over Y}\)  \widetilde{H}\({k^2 X\over n}\).
$$
Using Lemma 1, and a straightforward calculation, we obtain
$$
\sum\Sb n\le x\\ n\text{ odd} \endSb \frac{G_{k^2}(n)}{\sqrt{n}} 
= \frac{4}{\pi^2} x + O(\log x).
$$
Hence, by partial summation,
$$
M_k ={4\over\pi^2}\intl_0^\infty \Phi\(t \over Y\)
\,\widetilde{H}\({k^2 X\over t}\) 
{dt\over\sqrt{ t}} + O \biggl( \log Y \intl_{0}^{Y} 
\biggl| \biggl(\Phi\(t\over Y\) \widetilde{H} \( {k^2 X}\over 
{t}\) \frac{1}{\sqrt{t}}\biggr)^{\prime}\biggr|dt\biggr).
$$
Making the change of variable $t=yY$ and
using the bounds (4.1) (with $j=1$ or $2$),
we have 
$$
M_k = \frac{4}{\pi^2} \sqrt{Y}
\intl_0^{\infty} \Phi(y) \widetilde H\( k^2X \over {yY}
\) \frac{dy}{\sqrt{y}} +O\biggl(\frac{U\sqrt{Y}}{k^2 X}\log Y
\biggr).
$$

Combining this with our earlier expression for $M_0$ we
conclude that
$$
M = \frac{X\sqrt{Y}}{\pi^2} 
\intl_0^\infty {\Phi(y)\over \sqrt{ y}}
\sum_{k=-\infty}^\infty 
\,\widetilde{H}\({k^2X \over {yY}}\) 
dy + O\((X+U\sqrt{Y}) \log Y\). 
\eqno(5.1)
$$
Using (4.1) (with $j=1$) when $|k|\ge \sqrt{UY}/\sqrt{X}$, and (4.1) for 
smaller $k$ we deduce that 
$$
\sum_{k=-\infty}^{\infty} \widetilde{H}\({k^2X} \over { yY}\) 
= 1+\frac{yY}{2\pi X} 
\sum\Sb k=-\infty\\ k\neq 0\endSb^{\infty} \frac{1}{k^2}
\biggl(1-\cos \(\frac{2\pi k^2X}{yY}\) +
\sin\({2\pi k^2X \over yY}\)\biggr) +O\biggl(\frac{\sqrt{Y}}
{\sqrt{UX}}\biggr).
$$
Using this in (5.1) we obtain
$$
M= \frac{2X^{3\over 2}}{\pi^2} C\(\frac YX \) + O\biggl(\frac{X^{3\over 2}
+Y^{3\over 2}}{\sqrt{U}} + (X+U\sqrt{Y})\log Y\biggr). 
$$
We combine this with (3.1), (3.2) and (4.3), and choose $U=(XY)^{2\over 5} 
(X+Y)^{-{3 \over 5}}$ to obtain Theorem~1.

\head
6. Some graphs of $C(\al)$
\endhead

In this section we present graphs of $C(\alpha)$ and $C'(\alpha)$,
and we prove that $C'(\alpha)$ is differentiable 
at $\alpha\in \Q$ if and only if $\alpha=2p/q$ with $p$ and $q$
both odd.

\noindent\vbox{
\vskip \baselineskip
\noindent The following graphs show $C(\alpha)$,

\vskip .1in

\hskip .1in \epsffile{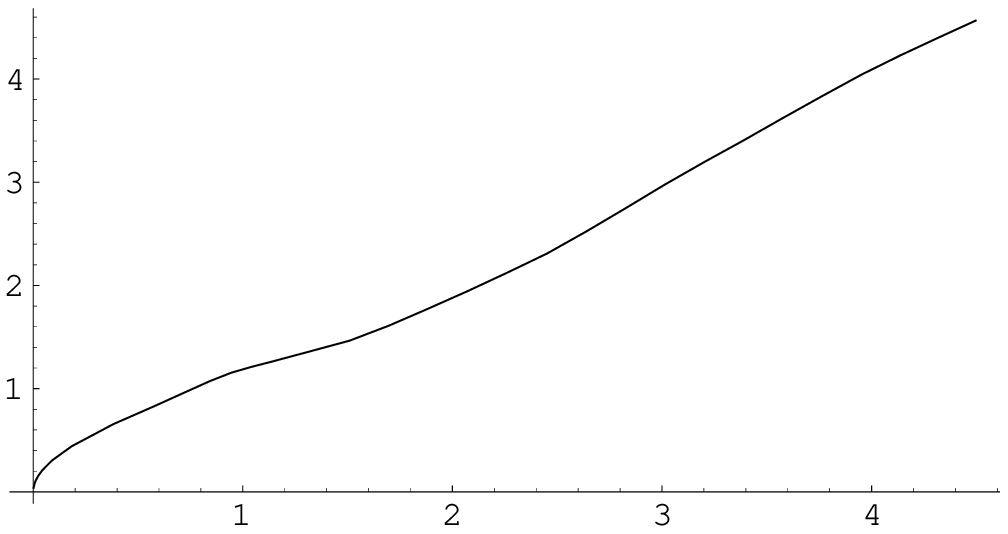}
}
\vskip .1in


\noindent\vbox{
\noindent and  $C'(\alpha)$:


\vskip .10in

\hskip .1in \epsffile{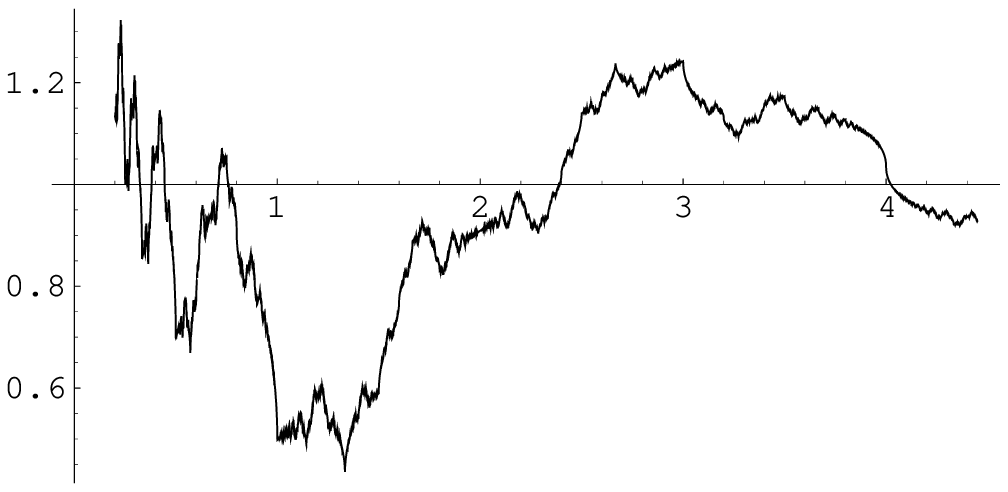} 
}
\vskip .1in

It certainly appears from the graph that $C'(\alpha)$ is
not everywhere differentiable.  

\proclaim{Proposition} $C'(\alpha)$ is differentiable
at $\alpha\in\Q$ if and only if $\alpha=2p/q$ with $p$
and $q$ both odd.
\endproclaim

The differentiability of $C'(\alpha)$ at $\alpha=2$ can be seen
in the above graph.

\demo{Proof}
Let 
$$
f(\alpha)=\frac 1{\pi} \sum_{k\not=0} {\sin(\pi k^2 \alpha)\over k^2}.
$$
By the second expression for $C(\alpha)$ in Theorem~1 we have
$$
(\alpha^{-{3\over 2}}C(\alpha))'=-\frac 12 \alpha^{-{3\over 2}}
-
\alpha^{-\frac 52 }
f\(\frac{\mathstrut\alpha}2 \),
\eqno{(6.1)}$$
the term--by--term differentiation being justified by the 
uniform absolute convergence of the resulting sum.
We will show that $f(\alpha)$ is
differentiable at ${p\over q}\in\Q$ if and only if
$p$ and $q$ are both odd.   

The following Lemma is an exercise in the Poisson summation
formula.

\proclaim{Lemma} Suppose $(p,q)=1$.
As $\alpha\to 0^\pm$ we have
$$
f\({p\over q}+\alpha\)=f\({p/q}\)
-\alpha
\pm
\sqrt{|2\alpha|}\({1\over 2q}\sum_{v\mod 2q}
\cos {\pi pv^2\over q}
\mp
\sin {\pi pv^2\over q}
\)
+O(\alpha^{3\over 2}).
$$
In particular, $f(\alpha)$ is differentiable at $\alpha={p\over q}$
if and only if $G(p/q)=0$, where
$$
G(p/q):=\sum_{v\mod 2q} e^{\pi i p v^2/q}.
$$
\endproclaim

The value of $G(p/q)$ is well known and can be
found in Chapter~2 of Davenport \cite{D}.  
We find that $G(p/q)=0$   if and only if $p$ and $q$ are both odd.
Combining this with formula (6.1) completes the
proof of the Proposition. 
Note that $C'(\alpha)$ is right-- or
left--differentiable at certain other $p/q$, such as the odd integers.
This can be determined by
considering the real and imaginary parts of~$G(p/q)$.
\enddemo

\Refs

\item{[D]} {\sl H. Davenport}, Multiplicative Number Theory 2nd ed,
GTM 74, Springer--Verlag, 1980.

\item{[Du]} {\sl J.J.~Duistermaat}, Selfsimilarity of Riemann's
nondifferentiable function, Nieuw archief voor wiskunde {\bf 9}
(1991), 303--337.

\item{[G]}  {\sl J.~Gerver}, The differentiability of the Riemann 
function at certain rational multiples of $\pi$, 
Amer. J. Math. {\bf 92} (1970), 33--55.

\item{[H]} {\sl G.H.~Hardy}, Weierstrass's non--differentiable function,
Trans. Amer. Math. Soc. {\bf 17} (1916), 322--323.

\item{[J]} {\sl M.~Jutila}, On the mean value of $L(1/2,\chi)$
for real characters, Analysis {\bf 1} (1981), 149--161.


\item{[S]} {\sl K.~Soundararajan},  Nonvanishing of quadratic
Dirichlet $L$--functions at $s=\frac 12$, preprint.

\item{[Sm]} {\sl A.~Smith}, Differentiability of Riemann's function,
Proc. Amer. Math. Soc. {\bf 34} (1972), 463--468.

\endRefs

\enddocument

\head
7. Details of nondifferentiability
\endhead 

[[I suggest not including this section in the paper, 
but I wanted to
write these things down so they can be checked.]]
[[This is probably an exercise in some book.  Can anyone
provide a reference?]]

First we use the Poisson summation formula a few times.
Suppose $\alpha>0$ and $0\le v < u$. 
$$\eqalign{
{\sum_k}' {\sin\pi (uk+v)^2\alpha\over (uk+v)^2}
&=-\delta_{0v}\pi \alpha
+\sum_j \widehat{h_{s,\alpha}}(j)\cr
&=-\delta_{0v}\pi \alpha
+ {\sqrt{|\alpha|}\over u}
\sum_j e(jv/u) \widehat{h_s}\({j\over u\sqrt{|\alpha|}}\)\cr
&=-\delta_{0v}\pi \alpha
+{\pi\sqrt{|2\alpha|}\over u} 
+O(u\alpha^{3\over 2})}$$
where
$$\leqalignno{
h_{s,\alpha}(x)&={\sin \pi (ux+v)^2\alpha\over (ux+v)^2}\cr
&&\hbox{and}\cr
h_s(x)&={\sin \pi x^2\over x^2}.
}$$
Note that $\widehat{h_s}(0)=\pi\sqrt{2}$ and we can (barely)
integrate by parts twice to show
$\widehat{h_s}(j)\ll j^{-2}$.  If $\alpha<0$ use the fact that
the above is an odd function.

Also
$$\eqalign{
{\sum_k}' {\cos\pi (uk+v)^2\alpha\over (uk+v)^2}
&={\sum_k}'{1\over (uk+v)^2} + 
{\sum_k}' {\cos\pi (uk+v)^2\alpha-1\over (uk+v)^2}\cr
&={\sum_k}'{1\over (uk+v)^2} + 
\sum_j \widehat{h_{c,\alpha}} (j)\cr
&={\sum_k}'{1\over (uk+v)^2}
+ {\sqrt{|\alpha|}\over u}
\sum_j e(jv/u) \widehat{h_c}\({j\over u\sqrt{|\alpha|}}\)\cr
&={\sum_k}'{1\over (uk+v)^2}
-{\pi\sqrt{|2\alpha|}\over u} 
+O(u\alpha^{3\over 2})}$$
where
$$\leqalignno{
h_{c,\alpha}(x)&={\cos \pi (ux+v)^2\alpha-1\over (ux+v)^2}\cr
&&\hbox{and}\cr
h_c(x)&={\cos \pi x^2-1\over x^2}.
}$$
Note that $\widehat{h_c}(0)=-\pi\sqrt{2}$ and we can (barely)
integrate by parts twice to show
$\widehat{h_c}(j)\ll j^{-2}$.
If $\alpha<0$ use the fact that
the function is even.

We have shown:
\proclaim{Lemma}  As $\alpha\to 0^\pm$ we have
$${\sum_k}' {\sin\pi (uk+v)^2\alpha\over (uk+v)^2}
=-\delta_{0v}\pi \alpha
\pm{\pi\sqrt{|2\alpha|}\over u} 
+O(u\alpha^{3\over 2})
$$
and
$$
{\sum_k}' {\cos\pi (uk+v)^2\alpha\over (uk+v)^2}
=A_{u,v}-{\pi\sqrt{|2\alpha|}\over u} 
+O(u\alpha^{3\over 2})
$$
for some constant $A_{u,v}$.
\endproclaim

\demo{Proof of the Lemma in Section 6}
Note that, no matter what the parity of $p$ or $q$,
both $\cos\pi k^2p/q$ and $\sin\pi k^2p/q$ only
depend on $k\mod 2q$.  
We have
$$\eqalign{
\sum_{k\not=0}
{\sin\pi k^2(\alpha+p/q)\over k^2} =&
\sum_{k\not=0}
{\sin\pi k^2\alpha\cos\pi k^2p/q\over k^2}+
{\cos\pi k^2\alpha\sin\pi k^2p/q\over k^2}\cr
=&\sum_{v\mod 2q} \cos (\pi p v^2/q)
{\sum_{k}}' {\sin\pi (2qk+v)^2\alpha\over (2qk+v)^2}
\cr
&+
\sum_{v\mod 2q} \sin (\pi p v^2/q)
{\sum_{k}}' {\cos\pi (2qk+v)^2\alpha\over (2qk+v)^2}.
}$$
Insert the above lemma to get the expression in Section~6.
We used $\sum \sin(\pi pv^2/q)A_{2q,v}=f(p/q)$,
which follows from the continuity of $f$.
\enddemo

In the last step we used the
following values of the Gauss sum:
$$
{1\over 2\sqrt{q}}\sum_{v\mod 2q} e^{\pi i {p\over q}v^2}=
\cases
0 & p,\ q\ odd\cr
(1+i)/\sqrt{2}& q\ even,\ p=\square\mod 2q\cr
-(1+i)/\sqrt{2}& q\ even,\ p\not=\square\mod 2q\cr
1&q\equiv 1\mod 4,\ p/2=\square\mod q\cr
-1&q\equiv 1\mod 4,\ p/2\not=\square\mod q\cr
i&q\equiv 3\mod 4,\ p/2=\square\mod q\cr
-i&q\equiv 3\mod 4,\ p/2\not=\square\mod q\cr
\endcases
$$
I think I checked all possibilities.

\head
8. Twisting by  $\ell$
\endhead

[[Again, I suggest not including this section in the paper, 
but I wanted to
write these things down so they can be checked.]]

Everything works the same as in the $\ell=1$ case up to
the point
$$\eqalign{
\Cal S_\ell (X):=&\sum_{m \ odd\ } \sum_{n\ odd} 
\js m{\ell n} H\({n/ X}\) \Phi\({m/ X}\) \cr
=& M(\ell)+R(\ell)
}$$
where        
$$\eqalign{                                                                  
M(\ell):=&{X\over 2\ell}\sum_{k=0}^\infty \sum_{n\ odd} 
H\({n/ X}\) 
{G_{k^2}(\ell n)\over n} 
\,\widetilde{\Phi}\({k^2X\over n\ell }\) 
\cr 
=&{X\over 2\ell}\sum_{k=0}^\infty M_k(\ell), \cr
}$$
say.  

\proclaim{Proposition} Suppose $\ell$ is odd and squarefree.
We have
$$
\sum_{n\le x\atop n\ odd} G_{k^2}(\ell n) \sim
\cases
\displaystyle
{4\over 3\pi^2}{\ell^{3\over 2}\over \sigma(\ell)} x^{3\over 2}
& if\ k=0 \cr
\displaystyle
{8\over 3\pi^2}{\ell^{3\over 2}\over \sigma(\ell)} x^{3\over 2}
& if\ k\not=0, \cr
\endcases
$$
where $\sigma(\ell)$ is the divisor sum function.
\endproclaim

By partial summation we obtain, for $k\ge 1$,
$$
M_k(\ell)\sim {4\over \pi^2}{\ell^{3\over 2}\over \sigma(\ell)}
X^{1\over 2}\intl_{-\infty}^\infty {H(t)\over \sqrt{t}}  
\widetilde{\Phi}\({k^2\over t\ell}\) dt.
$$
So,
$$\eqalign{
M(\ell)&\sim 
{X^{3\over 2}\over \pi^2}
{\ell^{1\over 2}\over \sigma(\ell)}
\intl_{-\infty}^\infty {H(t)\over \sqrt{t}} 
\sum_{k=-\infty}^\infty 
\widetilde{\Phi}\({k^2\over t\ell}\) dt
\cr
&={X^{3\over 2}\over \pi^2}
{1\over \sigma(\ell)}
\intl_{-\infty}^\infty {H(t/\ell)\over \sqrt{t}} 
\sum_{k=-\infty}^\infty 
\widetilde{\Phi}\({k^2\over t}\) dt
.
}$$

Applying Poisson summation to the sum over $k$ gives
$$
M(\ell)\sim {X^{3\over 2}\over \pi^2}
{1\over \sigma(\ell)}
\intl_{-\infty}^\infty H(t/\ell)
\sum_{k=-\infty}^\infty 
\intl_{-\infty}^\infty 
{\Phi(u)\over \sqrt{u}}\cos\!\({\pi k^2 t\over 2u}\) du\,dt
$$

The details in the proof of the Proposition are in
the following Lemma.  There are a number of cases to consider,
because of possible common factors of $k$ and $\ell$.

\proclaim{Lemma} Let 
$$
\Pi_{k,\ell}(s)=\sum_{j=0}^\infty {G_{k}(\ell p^j)\over
p^{js} }.
$$
If $p\not|\ell$ and $p\not|k$  then 
$$\Pi_{k^2,\ell}\(s\)=G_{k^2}(\ell)\(1+{\sqrt{\mathstrut p}\over p^s}\).
$$
If $p\not|\ell$ then 
$$\Pi_{k^2,\ell}\(s\)=
G_{k^2}(\ell)
1_{3\over 2}(s)
\(1+{1\over p}\) .
$$
If $p|\ell$ and $p\not|k$ then
$$\Pi_{k,\ell}(s)=G_{k}(\ell). 
$$
If $p|k$ and $\ell=p^2$ then
$$\Pi_{k^2,p^2}(s)= 
G_{k^2}(p^2)
1_{3\over 2}(s)
\, p^2 ,
$$
where $1_{3\over 2}(s)$ is an Euler factor which is\/  $1$ at $s={3\over 2}$.
\endproclaim

\demo{Proof}  In the first case we use the multiplicativity of
$G_k$ and the following cases of Lemma~2.3:  if $p^\lambda||k$
and $j\ge 1$ then
$$
G_{k^2}(p^j)=\cases
\varphi(p^j)& if\ j\le 2\lambda\ and\ j\ is\ even \cr
p^{2\lambda+{1\over 2}} & if\ j=2\lambda+1 \cr
0 & otherwise.\cr
\endcases
$$
So we have
$$\eqalign{
\sum_{j=0}^\infty {G_{k}(p^j)\over p^{js} } 
&= \sum_{j=0}^\lambda {\varphi(p^{2j})\over p^{2js}}
+ {p^{2\lambda+{1\over 2}}\over p^{(2\lambda+1)s}}
\cr
&=1+\(1-{1\over p}\)\sum_{j=1}^\lambda {p^{2j}\over p^{2js}}
+ {p^{2\lambda+{1\over 2}}\over p^{(2\lambda+1)s} }  
\cr 
&=1+\(1-{1\over p}\)
\( {1-p^{(2-2s)(\lambda+1)}\over 1-p^{2-2s}} -1\) 
+ {p^{2\lambda+{1\over 2}}\over p^{(2\lambda+1)s} } .   
}$$
So,
$$\eqalign{\Pi_{k^2,1}\({\textstyle{ 3\over 2}}\)&=
1+\(1-{1\over p}\) \({1-p^{-(\lambda+1)}\over 1-p^{-1}}                  
-1\)
+ {p^{2\lambda+{1\over 2}}\over p^{3\lambda+{3\over 2}} }\cr
&=1+{1\over p} .
}$$

In the last case, Lemma 2.3 says that if $p\not|k$ then
$G_{k}(p^j)=\sqrt{\mathstrut p}$ if $j=1$ and 0 if $j>1$.  
So at most
the $j=0$ term in the sum makes a contribution.
\enddemo

\bye